\documentclass[12pt]{article}
\usepackage{amssymb}
\usepackage{amsmath}
\usepackage{amsthm}
\usepackage{color}
\usepackage{comment}
\usepackage{geometry}		
\usepackage{graphicx}

\DeclareFontFamily{OT1}{pzc}{}
\DeclareFontShape{OT1}{pzc}{m}{it}{<-> s * [1.10] pzcmi7t}{}
\DeclareMathAlphabet{\mathpzc}{OT1}{pzc}{m}{it}

\let\originalleft\left
\let\originalright\right
\renewcommand{\left}{\mathopen{}\mathclose\bgroup\originalleft}
\renewcommand{\right}{\aftergroup\egroup\originalright}

\newtheorem{theorem}{Theorem}[section]

\newtheorem{corollary}[theorem]{Corollary}
\newtheorem{lemma}[theorem]{Lemma}
\newtheorem{proposition}[theorem]{Proposition}

\theoremstyle{definition}
\newtheorem{definition}{Definition}[section]

\theoremstyle{remark}

\begin{document}

\title{Scaling of the rotation number\\ for perturbations of rational rotations 
}
\author{
Paul Glendinning\\
Deparment of Mathematics,\\ University of Manchester,\\
Manchester M13 9PL, U.K.
}
\maketitle


\begin{abstract}
The parameter dependence of the rotation number in families of circle maps which are perturbations of rational rotations is described. We show that if, at a critical parameter value, the map is a (rigid) rotation $x\to x+\frac{p}{q}~({\rm mod}~1)$ with $p$ and $q$ coprime, then the rotation number is differentiable at that point provided a transversality condition holds, and hence that the rotation number scales linearly at this parameter. We provide an explicit and computable expression for the derivative in terms of the Fourier series of the map, and illustrate the results with the Arnold circle map and some modifications. Piecewise linear circle maps can also be treated using the same techniques. 
\end{abstract}

\maketitle

\begin{quotation}
Circle maps provide simple models of systems with two interacting frequencies. As such they provide insights into
problems in areas such as climate change \cite{Ghil2020} and biological cycles (e.g. the heart beat \cite{Glass1991} and sleep-wake cycles \cite{DGS2021}). There are typically two important parameters: the ratio between the two fundamental frequencies and the size of the nonlinearity. The dynamics of these maps explains phenomena such as mode locking, quasiperiodicity and chaos. Mode locking and quasiperiodicity can be described via the rotation number, a measure of the average speed of rotation around the circle. For invertible maps, rational rotation numbers imply periodic motion and mode locking, typically on intervals of parameters. Similarly, irrational rotation numbers correspond to quasiperiodic behaviour, with the motion typically dense in the circle. In applications the scaling of the rotation number is often important and its variation with a parameter often forms a `devil's staircase', with the steps formed by intervals of mode locked behaviour. In the mode locked interval the rotation number is constant, but at the end of a mode locked interval the rotation number typically scales like the square root of the distance from the end of the interval \cite{Baesens1991,MT1986,PM1980}, and is therefore not differentiable there. In this paper we consider perturbations from a rigid rational rotation, $x\to x+\frac{p}{q}~({\rm mod}~1)$. If a simple transversality condition holds we show that the rotation number is differentiable at this special point, and hence the local scaling is linear, in contrast to the standard case. This clarifies comments by Herman \cite{H1977} on the non-differentiability of the rotation number in the rational case described in the context of the irrational case, extending classical results of Brunovsky \cite{Brunovsky1974} to the case of rational rotation number. We also provide an explicit formula for the derivative and apply it to some standard (and less standard) examples.
\end{quotation}

\section{\label{sect:intro}Introduction}

The rotation number is a topological invariant of orientation preserving homeomorphisms of the circle. In a series of results starting with Poincar\'e in the 1880s and continued by mathematicians such as Denjoy in the first half of the twentieth century and Herman in the second half of the twentieth century, it was shown that many properties of such maps depend on their rotation number $\rho$. One interesting issue is the extent to which the rotation number determines the conjugacy class of the map. Two maps $f$ and $g$ are topologically conjugate if there exists a homeomorphism $h$ such that $f\circ h=h\circ f$; $h$ can be thought of as a change of coordinates. If $h$ is a diffeomorphism then the maps are differentiably conjugate to each other. A flavour of the results proved is contained in the following statements about orbits, i.e. solutions of the difference equation $x_{n+1}=f(x_n)$.
\begin{itemize}
\item[(i)] $\rho\in\mathbb{R}\backslash\mathbb{Q}$ if and only if $f$ has no periodic orbits (Poincar\'e);
\item[(ii)] if $\rho\in \mathbb{R}\backslash\mathbb{Q}$ and $f$ is a $C^2$ diffeomorphism then $f$ is topologically conjugate to the rotation $x+\rho~(\rm{mod}~1)$ (Denjoy); and
\item[(iii)] if $\rho\in \mathbb{R}\backslash\mathbb{Q}$ and satisfies a Diophantine condition, and $f$ is a $C^3$ diffeomorphism, then $f$ is differentiably conjugate to the rigid rotation $x+\rho~(\rm{mod}~1)$ (Herman).
\end{itemize}
See e.g. \cite{Khanin2003,dMvS} for historical details. This class of result has been extended to piecewise smooth circle maps (e.g. \cite{Khanin2003}) and to the bifurcations in families of maps (e.g. \cite{Boyland1986,Quas1992,RT1991}). In such families the rotation number typically forms a devil's staircase as a function of parameters with plateaus at rational values corresponding to intervals of parameters, the mode locked intervals, on which the (rational) rotation number is constant. 

In standard smooth examples such as the Arnold map, 
\begin{equation}
F(x,\alpha ,\beta )=x+\alpha+\beta\sin (2\pi x).
\label{eq:Arn00}
\end{equation} 
(see Sec.~\ref{sect:Arnold}), the rotation number has mode locked plateaus and is not differentiable at the end-points of these plateaus as shown in Figure~\ref{fig:Arn}a of Sec.~\ref{sect:Arnold}. On the other hand, Brunovksy \cite{Brunovsky1974} and Herman \cite{H1977} prove that for a smooth family passing through an irrational rigid rotation, $x+\rho$, $\rho\in \mathbb{R}\backslash\mathbb{Q}$, the rotation number is differentiable at this point. In general this is not the case for families passing through rational rigid rotations, as Herman points out. However, in a wide variety of smooth families through rigid rational rotations, including families through parameter space of the Arnold map, we show here that the rotation number \emph{is} differentiable at rigid rational rotations. The aim of this note is to explain this phenomenon and provide a transversality condition  that implies this differentiability, see Theorem~\ref{thm:main} below. Obviously there is a trivial case where the rigid rotation occurs in the interior of a mode locked interval and the rotation number is constant with derivative zero. This is less surprising and is treated briefly at the end of the proof in Sec.~\ref{sect:proofmain}.

Some notation is needed to describe the results, see e.g. \cite{Devaney1989,dMvS} for more detail. The circle, $\mathbb{T}^1$ can be identified with the interval $[0,1)$ where we may think of 1 being identified with 0. An orientation preserving circle homeomorphism $f:\mathbb{T}^1\to\mathbb{T}^1$ is described by a lift $F:\mathbb{R}\to\mathbb{R}$ which is a continuous increasing function such that $f(x~({\rm mod}~1))=F(x)~{\rm mod}~1$. Two lifts of $f$ differ by an integer, and $F(x+1)=F(x)+1$, so lifts defined on $[0,1)$ can then be extended to $\mathbb{R}$ using this relationship. The rotation number of a lift $F$ is defined by
\begin{equation}
\rho (F)=\lim_{n\to\infty}\frac{1}{n}\left(F^n(x)-x\right),
\label{eq:rotnodef}
\end{equation}
The limit converges uniformly to a value independent of $x$, and the rotation numbers of two lifts differ by an integer \cite{Devaney1989}. If $F$ is a lift of $f$ then $F^n$ is a lift of $f^n$ ($f$ composed with itself $n$ times) and $\rho (F^n)=n\rho (F)$ (this will be useful later in this paper). The \emph{rigid rotation} with rotation number $\omega$ is the map with lift
\begin{equation}
R_\omega (x) = x +\omega .
\label{eq:rigidlift}\end{equation}
The lift of any orientation preserving circle homeomorphism can be written in the form $F(x)=x+\phi (x)$ where $\phi$ is periodic with period one and $\frac{\phi (y)-\phi (x)}{y-x}>-1$ if $y>x$. It is sometimes useful to separate out the mean of $\phi$ over $[0,1]$ and write $\phi (x)=a+\psi (x)$ with $\int_{\mathbb{T}^1}\psi (x)\, dx =0$.    

We will consider families of circle maps parameterised by a real variable $\mu$. If the family is continuous in the $C^0$ topology then the rotation number as a function of $\mu$, $\rho (\mu )$ with the obvious abuse of notation, is continuous \cite{dMvS,RT1991}. Let $D^1(\mathbb{T}^1)$ denote the space of lifts of orientation-preserving circle diffeomorphisms  

Herman's result, 5.3 of \cite{H1977}, with part (b) due to  Brunovsky \cite{Brunovsky1974}, uses the notation $f={\rm Id}+\alpha+\phi$ for the circle homeomorphisms, and where a $C^1$ path of lifts is a path in one-dimensional parameter space $f(x,\mu )=x+\alpha +\phi (x,\mu)$ such that for each fixed $\mu$, $\phi (x,\mu )\in D^k(\mathbb{T}^1)$ for some $k\ge 0$ and the map $\mu \to \phi(x,\mu )$ is $C^1$ with derivative $\frac{\partial \phi }{\partial \mu}(x,\mu )$. Herman \cite{H1977} denotes the parameter by $t$ but we use $\mu$ to emphasise the link with bifurcation theory. With the obvious abuse of notation we will write $\rho (f(x ,\mu))=\rho(\mu )$ defined via a $C^1$ path in the lifts of $f$.

\begin{theorem}\cite{H1977} Suppose $\alpha\in \mathbb{R}\backslash\mathbb{Q}$ then
\begin{itemize} 
\item[(a)] $\rho:D^1(\mathbb{T}^1)\to \mathbb{R}$ is differentiable at $R_\alpha$ and the derivative map is
\[
\phi \in C^1(\mathbb{T}^1)\to \int_{\mathbb{T}^1}\phi (x)\,dx\ \in\mathbb{R}
\]
($\rho$ is not differentiable at $R_{p/q}$, $p/q\in\mathbb{Q}$).
\item[(b)]  If $\mu\in [a,b]\to f_\mu(x)=x+\alpha +\phi_\mu(x)\in D^0(\mathbb{T}^1)$ is a $C^1$ path and $\phi_0(x)\equiv 0$, then $\mu\to \rho (\mu)$ is differentiable at $\mu=0$ and we have:
\[
\frac{d\rho}{d \mu}(\mu)\big|_{\mu =0}=\int_0^1\frac{\partial f}{\partial \mu}(x,\mu )\big|_{\mu=0}\, dx.
\]
\end{itemize}
\label{thm:H1977}
\end{theorem}

The focus in this note is on the remark in parentheses at the end of part (a) of the Theorem. We will work with a slightly simpler (but less general)  definition of maps which makes the derivative of the maps with respect to $\mu$ more transparent.

\begin{definition}We will say that a family of lifts is in class $\mathcal{R}_{p,q}$ in a neighbourhood of $\mu =0$ if there exists $\mu_0>0$ such that for all $|\mu |<\mu_0$
\begin{equation}
F(x,\mu )=x+\frac{p}{q}+\mu (a+\psi (x))+\mu^2g(x,\mu ),
\label{eq:Fclassq}\end{equation}
where $p$ and $q$ are coprime, $\psi(x)$ is a $C^1$ periodic function with period one with $\int_0^1\psi (x)\,dx=0$, and $g(x,\mu )$ is continuous in $\mu$, and for each $\mu\in (-\mu_0,\mu_0)$,  $g(x,\mu )$ is a $C^1$ periodic function of $x$ with period one. 
\label{def:Rpq}
\end{definition}

Herman's result is clearly applicable to families of diffeomorphisms of the circle with lifts in $\mathcal{R}_{p,q}$ in a neighbourhood of $\mu =0$ and since 
\begin{equation}
\frac{\partial F}{\partial \mu}(x,0)=a+\psi(x),
\label{eq:dFmu}\end{equation}
Brunovsky's result (part (b) of Theorem~\ref{thm:H1977}) can be rephrased as the following corollary. Here and below we will denote derivatives with respect to $\mu$ by primes.

\begin{corollary}\cite{Brunovsky1974} Suppose $\alpha\in \mathbb{R}\backslash\mathbb{Q}$ and $F(x,\mu )$ is a family of lifts in $\mathcal{R}_{p,q}$ then
$\rho (\mu)$ is differentiable at $\mu =0$ and
\begin{equation}
\rho'(0)=a.
\label{eq:brun0}\end{equation}
\label{cor:Brun}
\end{corollary}

The main result of this note is that there is an almost precisely equivalent formulation in the rational case provided a transversality condition holds.

\begin{theorem}Suppose $F(x,\mu )$ is a family of lifts in  $\mathcal{R}_{p,q}$  satisfying \eqref{eq:Fclassq} with $a>0$. Denote the Fourier coefficients of $\psi (x)$ by $(c_m)_{-\infty}^\infty$,  $c_0=0$, and define 
\begin{equation}
 \Psi (x)=\sum_{n\in\mathbb{Z}\backslash\{0\}}c_{nq}e^{2\pi i nqx}.
\label{eq:dFqmu}\end{equation} 
If 
\begin{equation}
\min_{x\in [0,1]}(a+\Psi(x))>0 ,
\label{eq:genericity}\end{equation} 
then $\rho (\mu )$ is differentiable at $\mu =0$ and    
\begin{equation}
\rho'(0)=T_0^{-1}, \quad T_0=\int_0^1\frac{1}{a+\Psi (x)}dx.
\label{eq:main}\end{equation}
If $\min (a+\Psi(x))<0<\max (a+\Psi(x))$ then there exists $\mu_1\in (0,\mu_0)$ such that $\rho (\mu )=\frac{p}{q}$ for all $|\mu|<\mu_1$ and so $\rho$ is differentiable at $\mu =0$ with $\rho'(0)=0$.
\label{thm:main}
\end{theorem}
   
If $\max (a+\Psi(x))<0$ then the conclusions hold after reversing the sign of $\mu$.

Theorem~\ref{thm:main} provides the counterpoint to Brunovsky's result \cite{Brunovsky1974}, part (b) of Theorem~\ref{thm:H1977} for rational rotation numbers. Note that if $c_{rq}=0$ for all $r\in\mathbb{Z}$ then \eqref{eq:main} becomes $\rho'(0)=a$, as in the irrational rotation number case, \eqref{eq:brun0}. Parkhe \cite{Parkhe2013}
proves a very similar result: that if the rotation number is strictly monotonic at $\mu =0$ then the rotation number is 
differentiable there. The differences between the approach here and that of \cite{Parkhe2013} are described in more detail at the end of this section. 

A sufficient condition for \eqref{eq:genericity} can be given in terms of properties of $\psi$ and $a$ without needing to find $\Psi$.

\begin{corollary}Suppose $F(x,\mu )$ is a family of lifts in  $\mathcal{R}_{p,q}$  satisfying \eqref{eq:Fclassq} with $a>0$. If $\min_{x\in [0,1]}(a+\psi (x) )>0$ then $\rho (\mu )$ is differentiable at $\mu=0$ with derivative given by \eqref{eq:main}.
\label{cor:qpsi}\end{corollary}  

We will illustrate our results with the Arnold sine map with lift \eqref{eq:Arn00}.
In particular we will prove two scaling results, the first corresponding to the case in which all coefficients of the Fourier series of $\psi$ of order $nq$ vanish, and the second a case in which they do not ($q=1$). 

\begin{theorem}Consider a $C^1$ path $(\alpha (\mu ),\beta (\mu ))$ through parameter space of \eqref{eq:Arn00} with $(\alpha (0),\beta (0))=(\frac{p}{q},0)$, $p,q$ coprime, $0\le p<q$. If $0<|\alpha '(0)|<\infty$ then $\rho (\mu )$ is differentiable at $\mu =0$ and if $q\ge 2$ then $\rho'(0)=\alpha '(0)$. If $(\alpha (0),\beta (0))=(0,0)$ and $|\alpha '(0)|>|\beta '(0)|$ then $\rho (\mu )$ is differentiable at $\mu =0$ with
\begin{equation*}
\rho'(0)={\rm sign}(\alpha '(0))\sqrt{|\alpha '(0)|^2-|\beta '(0)|^2}.
\end{equation*}
\label{thm:arnall}
\end{theorem}

See Corollary~\ref{cor:Arnold} and Corollary~\ref{cor:0diffA}.

In Sec.~\ref{sect:FS} we show that near a rational rotation with rotation number $\frac{p}{q}$ the $q^{th}$ iterate of the lift takes a simple form and is equivalent to the lift of a map with rotation number zero, i.e. lifts in  $\mathcal{R}_{0,1}$ and establish Corollary~\ref{cor:qpsi}. In Sec.~\ref{sect:Arnold} an immediate corollary is obtained that implies the $q\ge 2$ cases of Theorem~\ref{thm:arnall}.  Sec.~\ref{sect:int} considers the case $q=1$ using the Euler method for differential equations to bound solutions of the corresponding difference equation. This section  includes a corollary proving the $q=1$ case of Theorem~\ref{thm:arnall}. Then Sec.~\ref{sect:proofmain} brings together the results of Sec.~\ref{sect:FS} and Sec.~\ref{sect:int} to complete the proof of Theorem~\ref{thm:main}. In Sec.~\ref{sect:morearnold} a generalisation of the Arnold map with higher order harmonics is analysed showing how \eqref{eq:main} can be used to derive analytic forms of the scaling in more general situations. Sec.~\ref{sect:PWL} shows how the methods can be extended to the study of piecewise linear maps, providing a simple proof of a theorem in \cite{GMM2025}.

As noted earlier, Theorem~\ref{thm:main} has overlap with the main result of Parkhe \cite{Parkhe2013}. Parkhe shows that the rotation number is differentiable in a slightly more general setting, in that the transversality condition \eqref{eq:genericity} is replaced by the condition that the rotation number is strictly monotonic at $\mu =0$, and the error terms in \eqref{eq:Fclassq} are not made explicit, allowing the results to hold in a more general context. On the other hand, the derivative of the rotation number is left as
\begin{equation}
\rho'(0)=\left(q\int_0^1\frac{1}{\frac{\partial F^q}{\partial \mu}(x,0)}dx\right)^{-1},
\label{eq:Parkherot}\end{equation}  
which is much harder to work with than \eqref{eq:main} which provides an explicit integral in terms of the Fourier
series of $F$. There is inevitably some overlap in the method of proof: taking the $q^{th}$ iterate and computing the 
passage time though one unit. However, our use of the Euler approximation to derive the time of passage used in Sec.~\ref{sect:int} here is both intuitively appealing and simplifies the arguments. It may also have application to other problems. The examples of Secs.~\ref{sect:morearnold} and ~\ref{sect:PWL}
illustrate the power of the representation of the derivative via \eqref{eq:main}. 

  
\section{Reduction to the $q=1$ case}\label{sect:FS}
Recall the form of the families of lifts we are using, \eqref{eq:Fclassq}, and note that since $\psi$ has mean zero there is no constant term in its Fourier series; this role is played by the constant $a$. 

\begin{proposition}Suppose that $F(x,\mu )$ is a family of lifts in class  $\mathcal{R}_{p,q}$ given by \eqref{eq:Fclassq} with
\begin{equation}
\psi (x)=\sum_{m\in\mathbb{Z}\backslash\{0\}} c_m e^{2\pi i mx}.
\label{eq:FSpsi}
\end{equation}
Let
\begin{equation}
\Psi (x)=\sum_{n\in\mathbb{Z}\backslash\{0\}} c_{nq} e^{2\pi i nqx}.
\label{eq:FSPsi}
\end{equation}
There exists $\mu_0'$ such that if $|\mu |<\mu_0'$ then
\begin{equation}
F^q(x,\mu )=x+p+q\mu\left(a+\Psi(x)\right) +\mu^2\hat g (x,\mu ),
\label{eq:FqFS}
\end{equation}
where $\hat g$ is $C^1$, $\hat g (x+1,\mu)=\hat g(x,\mu )$ and $\hat g$ is continuous in $\mu$.
\label{prop:FS}
\end{proposition}

\emph{Proof:}

Since $F$ in class $\mathcal{R}_{p,q}$ (Definition~\ref{def:Rpq}) if $1\le k\le q$ then
\begin{equation}
\begin{array}{rl}
F^k(x,\mu)=&x+k\frac{p}{q}+ka\mu +\mu\sum_{r=0}^{k-1}\psi (F^r(x,\mu ),\mu)\\
&+\mu^2\sum_0^{k-1}g(F^r(x,\mu ),\mu) .
\end{array}
\label{eq:Fk}
\end{equation}
In particular, 
\begin{equation}
F^k(x,\mu)=x+k\frac{p}{q}+\mu \hat R_k(x,\mu ),
\label{eq:Fksimple}
\end{equation}
where the error term $\hat R_k$ is continuous in $\mu$ and
$C^1$ and periodic with period one in $x$.  Thus  from \eqref{eq:Fk} with $k=q$ and \eqref{eq:Fksimple} for the intermediate evaluations
\begin{equation}
\begin{array}{rl}
F^q(x,\mu )=&x+q\frac{p}{q}+qa\mu +\mu\sum_{r=0}^{q-1}\psi (x+\frac{rp}{q})\\
& +\mu^2 g_q(x,\mu ),
\end{array}
\label{eq:Fq}\end{equation}
where the error term $g_q$ is continuous in $\mu$ and periodic with period one in $x$ (it includes the corrections obtained by substituting \eqref{eq:Fksimple} with $k=r$ into the $\mu$ term of \eqref{eq:Fk} with $k=q$). 

From the Fourier series representation of $\psi$ \eqref{eq:Fq} can be rewritten as
\begin{equation}\begin{array}{rl}
F^q(x,\mu)= & x+p+qa\mu +\mu\sum_{r=0}^{q-1} \left(\sum_{m\in\mathbb{Z}\backslash\{0\}} c_me^{2\pi i m(x+\frac{rp}{q})}\right)\\ & +\mu^2g_q(x,\mu ).\end{array}
\label{eq:Fqfourier}
\end{equation}

Since $\psi$ is $C^1$ it has an absolutely convergent Fourier series, see e.g. \cite[Chapter I.6]{Katznelson2004} and hence we may reverse the order of the summation and find that
\[
\sum_{r=0}^{q-1} \left(\sum_{m\in\mathbb{Z}\backslash\{0\}} c_me^{2\pi i m(x+\frac{rp}{q})}\right)=
\sum_{m\in\mathbb{Z}\backslash\{0\}} c_me^{2\pi i mx}\left(\sum_{r=0}^{q-1}e^{2\pi im\frac{pr}{q}}\right) .
\]
If $m\ne nq$ for some $n\in\mathbb{Z}$ then the finite sum is zero, whilst if $m=nq$ then the finite sum equals $q$. Hence substituting back into \eqref{eq:Fqfourier} gives \eqref{eq:FqFS}.
\newline\rightline{$\square$}

To prove Corollary~\ref{cor:qpsi} we need to consider how some properties of $\psi$ imply results for $F^q$.

\begin{lemma}Let $F$, $q$, $\psi$ and $\Psi$ be as in Proposition~\ref{prop:FS}. Then $\min (a+\psi (x))>0$ implies that $\min (a+\Psi (x))>0$.
\label{lem:psitoPsi}\end{lemma}

\emph{Proof:}

If $\Psi \equiv 0$ then $\min(a+\Psi (x))=a\ge  \min(a+\psi(x) )>0$ and the lemma is trivially true. (Recall that $\psi$ has mean zero.)

If $\Psi$ is not identically zero then note that if $G_1$ and $G_2$ are lifts of circle homeomorphisms $g_1$ and $g_2$, and $G_1(x)<G_2(x)$, then $\rho (G_1)\le \rho (G_2)$. Suppose that $a-\psi (x)>m_1>0$, then there exists $B>0$ and $\mu_0>0$ such that if $0<\mu <\mu_0$ then $F(x,\mu )>x+\frac{p}{q}+m_1\mu+B\mu^2$ and so $\rho (F )>\frac{p}{q}$.

Now suppose that $\min (a+\Psi (x))<0$ and seek a contradiction (note that the inequality is strict since if $\min (a+\Psi (x))=0$ then we are in the $\Psi \equiv 0$ case). Since $\Psi$ has mean zero $\min (a+\Psi (x))<0$ implies that $\max (a+\Psi(x))>0$ and so provided $\mu_0$ is sufficiently small \eqref{eq:FqFS} implies that there exists $x_1$ and $x_2$ such that $F^q(x_1)-x_1<p<F^q(x_2)-x_2$. By the intermediate value theorem there exists $y$ such that $F^q(y)-y=p$ so $y$ (modulo one) is periodic and $\rho (F^q)=p$ or $\rho (F)=\frac{p}{q}$. This contradicts the statement that  $\rho (F )>\frac{p}{q}$ if $\min (a+\psi (x))>0$, hence $\min (a+\Psi (x))>0$.

If $\mu <0$ then inequalities are reversed but the same argument holds. 
\newline\rightline{$\square$}

Clearly Lemma~\ref{lem:psitoPsi} and Theorem~\ref{thm:main} imply Corollary~\ref{cor:qpsi}.


\section{The Arnold circle map}\label{sect:Arnold}
\begin{figure}[t!]
\centering
\includegraphics[width=4cm]{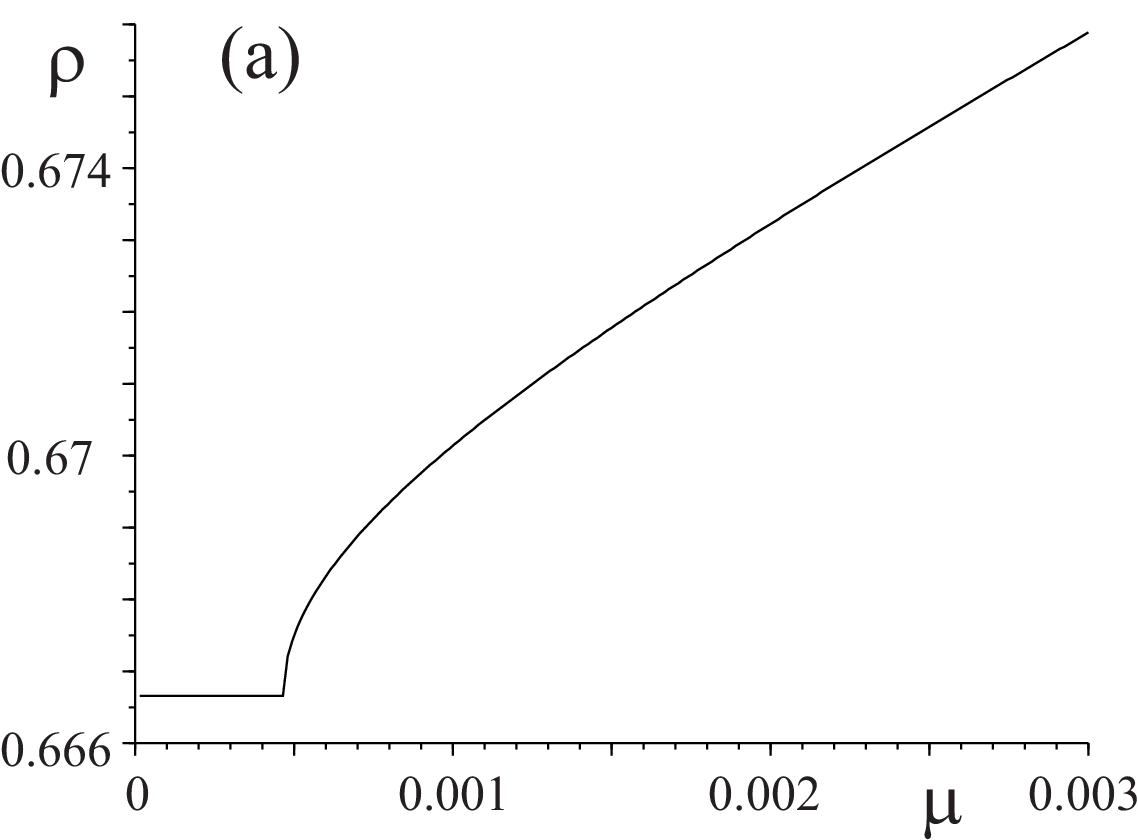}
\hspace{0.5cm}
\includegraphics[width=3.9cm]{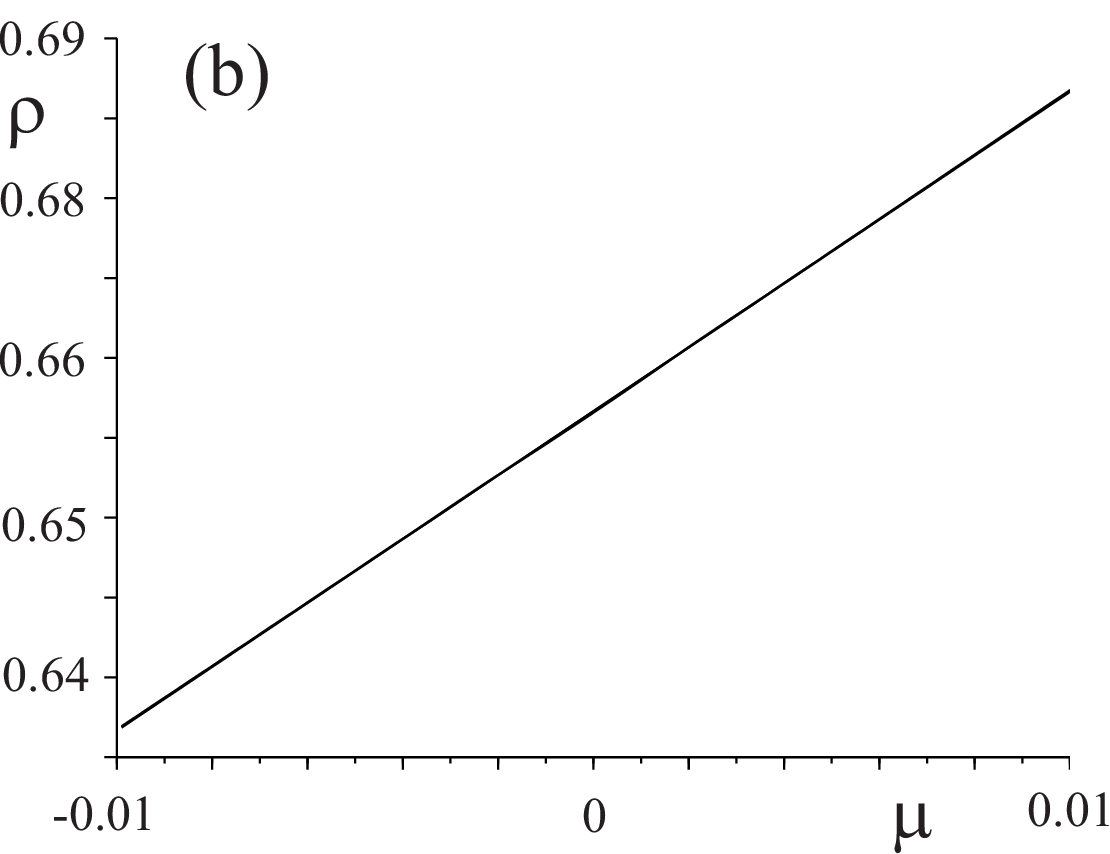}
\caption{
Rotation number $\rho$ of \eqref{eq:Arn00} against parameter $\mu$. Rotation numbers were computed using $10^5$ iterations computed in Scilab \cite{Scilab}. In both cases 200 data points are computed at equally spaced values of $\mu$. (a) $(\alpha (\mu),\beta (\mu))=(0.6615+2\mu, 0.1+\mu)$, $0<\mu <0.0003$, showing the typical non-differentiable transition at the end of an Arnold tongue. (b) Two curves superposed: $(\alpha (\mu),\beta (\mu))=(\frac{2}{3}+2\mu, \mu)$ and  $(\alpha (\mu),\beta (\mu))=(\frac{2}{3}+2\mu, 0.1\mu)$, $-0.01<\mu <0.01$. Both curves have slopes very close to 2 as expected from the theory (see main text).
}
\label{fig:Arn}
\end{figure}

The two-parameter family of circle maps with lifts \eqref{eq:Arn00}
is probably the most studied map of the circle in the literature. It is a homeomorphism if $|\beta |<\frac{1}{2\pi}$ and it is the rotation $R_\alpha$ if $\beta=0$, so we will be interested in one-parameter families through $(\alpha ,\beta )=(\frac{p}{q},0)$. Figure~\ref{fig:Arn}a shows the graph of the rotation number as a function of $\mu$ on a path that has $\beta >0$. This shows the nonlinear (square root) behaviour of the rotation number at the end of an Arnold tongue away from the line $\beta =0$ \cite[Sec. 4.4.1 and Appendix C.2]{Baesens1991}, and it is this general case that we want to contrast with the cases covered here. In Figure~\ref{fig:Arn}a $\beta$ is fixed at $\beta =0.1$. The end of the Arnold tongue with rotation number $\rho =\frac{2}{3}$ is close to $\alpha =\frac{2}{3}$ as $\beta$ is small, and the parameter window of Figure~\ref{fig:Arn}a starts at $\alpha =0.6615$ so that part of the plateau with rotation number $\frac{2}{3}$ is clearly visible. The non-differentiable transition at the end of the tongue is also clear, but for small values of $\beta$ the mode locked regions can be very thin, so the devil's staircase is not apparent at this resolution.   

The following is a corollary of the Fourier series description of Proposition~\ref{prop:FS}.

\begin{corollary}Consider a $C^1$ path parameterized by $\mu$ through $(\alpha ,\beta )=(\frac{p}{q},0)$ at $\mu=0$. If $q\ge 2$ and
\[
0<|\alpha'(0)|<\infty ,
\]
then the rotation number $\rho (\mu)$ is differentiable at $\mu =0$ and $\rho'(0)=\alpha '(0)$.
\label{cor:Arnold}
\end{corollary}

\emph{Proof:}
In the form \eqref{eq:Fclassq}, $\psi (x)=\beta \sin (2\pi x)$ and so all the Fourier coefficients of the Fourier series of $\psi$ with order greater than one are zero. Hence by Proposition~\ref{prop:FS} if $q\ge 2$ then 
\[
F^q(x,a,b )=x+q\alpha (0)+q\alpha'(0)\mu +O(\mu^2).
\]
Since $q\alpha (0)=p$, $F^q(x,a,b )=x+p+q\alpha'(0)\mu +O(\mu^2) $ and an elementary calculation with the $O(\mu^2)$ terms bounded by constants times $\mu^2$ in a small neighbourhood of $\mu =0$) implies that $\rho (F^q)=p+q\alpha '(0)\mu +O(\mu^2)$ and hence $\rho (F)=\frac{p}{q}+\alpha '(0)\mu +O(\mu^2)$. 
\newline\rightline{$\square$}

To illustrate Corollary~\ref{cor:Arnold} we will choose paths through $(\frac{2}{3},0)$ of the form
\begin{equation}
\alpha =\frac{2}{3}+v_1\mu , \quad \beta =v_2\mu ,  \quad |\mu |\ll 1 .
\label{eq:Apath}\end{equation}  
For such families 
\[
\frac{\partial F}{\partial \mu}=v_1+v_2\sin (2\pi x) ,
\]
so by Corollary~\ref{cor:Arnold}, the derivative of the rotation number at $\mu=0$ is $v_1$, independent of $v_2$. Figure~\ref{fig:Arn}b shows the curve of rotation number as a function of $\mu$ with $(v_1,v_2)$ equal to $(2,0.1)$ and $(2,1)$. As expected the two curves appear tangential at the intersection point, indeed, they are almost indistinguishable at this resolution. The linear regression routine, {\tt reglin} of Scilab \cite{Scilab} gives the slopes as 2.001 and 2.000 respectively (again as expected since $v_1=2$, and the result is independent of $v_2$) with a standard deviation of less than $3\times 10^{-5}$.


\section{Integer rotation number: the case $q=1$}\label{sect:int}
To prove the main result it is necessary to bound the rotation number as a function of the parameter. From Sec.~\ref{sect:FS}, Proposition~\ref{prop:FS}, the $q^{th}$ iterate of the map is a translation by an integer of the case $p=0$, $q=1$. Therefore we need to consider this case in more detail.

 One of the bounding techniques uses the comparison of a difference equation with a corresponding ordinary differential equation using Euler's method. Let us first recall the basic global error bounds for Euler's method for future reference.

Euler's method approximates solutions of the differential equation $\dot X=V(X)$ at times $nh$ by solutions of the difference equation $x_{n+1}=x_n+hV(x_n)$. If $L=\max |V'(x)|$, where the prime denotes the derivative with respect to $x$, and $M=\max |V(x)V'(x)|$ on the regions on which solutions are defined then (see almost any introductory textbook on numerical methods and differential equations, e.g. \cite[Chapter 1.13]{Braun1993})
\begin{equation}
|X(nh)-x_n|\le k_0h, \quad k_0=\frac{M}{2L}\left(e^{nhL}-1\right) .
\label{eq:Eulersmooth}
\end{equation}
  
The second elementary result needed relates the rotation number of the lift of a circle map to the iterates of a point (which can be taken to be $x=0$ without loss of generality).

\begin{lemma}Let $F:\mathbb{R}\to \mathbb{R}$ be the lift of an orientation preserving circle map. If there exists $n>0$ such that $F^n(0)\le 1\le F^{n+1}(0)$ then $\frac{1}{n+1}\le \rho(F)\le\frac{1}{n}$.
\label{lem:bddrot}
\end{lemma}

The proof uses standard techniques, see e.g. \cite{Devaney1989}, and is omitted.

From Sec.~\ref{sect:FS} the $q^{th}$ iterate of the maps $F$ we consider have rotation numbers close to zero. These maps are described in the next lemma.

\begin{proposition}Suppose $G$ is in class $\mathcal{R}_{0,1}$, i.e. 
\begin{equation}\label{eq:bigG}
G(x,\mu)=x+\mu (a + \psi(x))+\mu^2g(x,\mu ), 
\end{equation}
where $\psi$ and $g$ are periodic with period one and $C^1$ in $x$, and $g$ is continuous in $\mu$.  If $\min (a+\psi(x))>0$ then there exists $\mu_0>0$ and $k_1>0$ such that for all $0<\mu <\mu_0$
\begin{equation}
|\rho (G)-T_0^{-1}\mu |<k_1\mu^2, \quad T_0=\int_0^1\frac{1}{a+\psi(x)}dx.
\label{eq:rhoclose}
\end{equation}
\label{prop:rhoclose}
\end{proposition}

\emph{Proof:}

Let $m=\min (a+\psi (x))>0$ and $M=\max (a+\psi (x))$. Since $g$ is bounded, $|g|<B$ say. So we may choose $\mu_0>0$ such that $a+\psi (x)+\mu g(x,\mu ) \ge m-\mu_0B>0$ for all $x\in\mathbb{R}$ and all $\mu\in [0,\mu_0]$.  

Equation \eqref{eq:bigG} is the Euler iteration for solutions to the ordinary differential equation 
\begin{equation}
\frac{dX}{dt}=a + \psi(X)+\mu g(X,\mu ).
\label{eq:ode}\end{equation}
Since $\psi$ and $g$ are $C^1$ and $\epsilon >0$ solutions exist and are unique by Picard's Theorem (e.g. Ref.~\cite{Braun1993}) and the time taken to move from $X=0$ to $X=1$ is $T(\mu )$ where
\begin{equation} 
T(\mu )=\int_0^1\frac{1}{a + \psi(X)+\mu g(X,\mu )}dX.
\label{eq:Tmu}\end{equation}
The assumption $m-\mu_0B>0$ implies that this integral is finite.
Elementary manipulation of \eqref{eq:Tmu} with $\mu_0<\frac{1}{2}$ and using $1-x<(1-x)^{-1}<1+2x$ ($0\le x<\frac{1}{2}$) and the approximation $T_0\le \frac{1}{m}$ implies that 
\begin{equation}
|T(\mu )-T_0|<k_2|\mu |B,
\label{eq:T0Tmudiff}
\end{equation}
where $k_2=\min\{m^{-2},2m^{-1}M^{-1}\}$. 

The time $T(\mu )$ is derived from the ordinary differential equation \eqref{eq:ode}, but to translate this back into solutions of the difference equation \eqref{eq:bigG} it is necessary to take into account the global error of the Euler method, \eqref{eq:Eulersmooth}. 
 
The number of time steps in the Euler approximation of solutions to \eqref{eq:ode} is essentially $T(\mu)/\mu$, but generally this is not an integer so define $n(\mu )\in \mathbb{N}$ such that 
\begin{equation}
n(\mu)\mu \le T(\mu )<[(n(\mu)+1]\mu 
\label{eq:nmuTmuT}\end{equation}
with
\begin{equation}    
X(n(\mu)\mu)\le 1<X(n(\mu )+1)\mu).
\label{eq:nmuTmuX}\end{equation}
However, at $n(\mu)$ the solution of the difference equation \eqref{eq:bigG} from the solution is bounded by the Euler error \eqref{eq:Eulersmooth}. Define  $Q=M+B\mu_0$ and $P=m-B\mu_0 >0$ and note that $T(\mu )\le \int_0^1 P^{-1}dX=P^{-1}$. Then the Euler error is 
\[
|X(n(\mu )\mu )-x_{n(\mu)}|< \frac{1}{2}Q\mu(e^{n(\mu)\mu L}+1),
\]
where $L=\max(|\phi'(x)|)$. Since $n(\mu )\mu \le T(\mu)<P^{-1}$ this implies that
\begin{equation}
|X(n(\mu )\mu )-x_n|< k_3\mu, \quad k_3=\frac{1}{2}Q(e^{LP}+1).
\label{eq:Xx}\end{equation}
Also, from the definition of $Q$ and $P$, for all $r\ge 0$
\[
x_r+P\mu <x_{r+1}<x_r+Q\mu .
\]
Now choose  $k_4>0$ and less than the maximum of the right hand side of \eqref{eq:ode} with $\mu <\mu_0$.  Then by \eqref{eq:nmuTmuX}, $|X(n\mu )-1|<k_4\mu$, so there exists $J\le (k_3+k_4)/P$ independent of $\mu$ and $j(\mu)$ with $|j(\mu)|<J$ such that
\[
x_{n(\mu )+j(\mu)}\le 1<x_{n(\mu )+j(\mu)+1},
\]
and  $0< x_{n(\mu )+j(\mu)+1}-x_{n(\mu )+j(\mu)}<Q\mu $. Lemma~\ref{lem:bddrot} now implies that
\begin{equation}
\frac{1}{n(\mu )+j(\mu)+1}\le \rho (G)\le\frac{1}{n(\mu )+j(\mu)}.
\label{eq:rhoapprox}
\end{equation} 
To complete the proof we will use
\begin{equation}\begin{array}{rl}
\big|\rho(G)&-T_0^{-1}\mu \big|\le  \big|\rho(G)-\frac{1}{n+j+1}\big|+\ \big|\frac{1}{n+j+1}-\frac{1}{n}\big|\\
& \\
&+\ \big|\frac{1}{n}-T(\mu )^{-1}\mu\big|+\big|T(\mu )^{-1}\mu-T_0^{-1}\mu \big|.
\end{array}
\label{eq:bigineq}\end{equation}
(No doubt there are other ways to treat this but this is convenient.) The first step is to bound $n(\mu )$ in terms of $\mu$. From \eqref{eq:T0Tmudiff} and \eqref{eq:nmuTmuT}, if $\mu_0<\frac{T_0}{2(k_2+1)}$ then 
\begin{equation}
\frac{1}{n}<\frac{2\mu}{T_0}.
\label{eq:crudemun}
\end{equation}
(Note this is not a very accurate bound but is sufficient for the error estimates below.)
Each of the four terms on the right hand side of \eqref{eq:bigineq} can now be bounded using elementary manipulations which we indicate here without providing full details.

From \eqref{eq:rhoapprox} with $|j(\mu ||\le J$ and \eqref{eq:crudemun} we can bound the first term by 
\[
\frac{1}{n(\mu) -J}-\frac{1}{n(\mu )+J+1}\le \frac{8(2J+1)\mu^2}{T_0^2},
\]
provided $n>\sqrt{2}J$ or $\mu_0 < T_0/(2\sqrt{2}J)$.

The second term is straightforward: it is less than
\[
\frac{1}{n}-\frac{1}{n+J+1}\le \frac{4(J+1)\mu^2}{T_0^2}.
\]
To bound the third term use the reciprocal of \eqref{eq:nmuTmuT} to obtain an upper bound of $\frac{4\mu^2}{T_0^2}$.  Finally, multiplying out
\[
\frac{\mu}{T_0}-\frac{\mu}{T(\mu )}= \frac{(T(\mu )-T_0)\mu}{T_0T(\mu )}
\]
and using \eqref{eq:T0Tmudiff}, the fourth term is bounded by $\frac{2k_2\mu^2}{T_0^2}$  provided $\mu_0<\frac{1}{2k_2}$.

Putting these together, provided $\mu_0$ is sufficiently small so that each of the required inequalities to validate the bounds obtained above holds, we obtain \eqref{eq:rhoclose} with
\begin{equation}
k_1=T_0^{-2}(20J+16+2k_1).
\label{eq:kappa1}\end{equation}
\newline\rightline{$\square$}

This result has two useful corollaries.

\begin{corollary}If $G$ and $\hat G$ are in class $\mathcal{R}_{0,1}$ and 
\[
\frac{\partial G}{\partial \mu}(x,0)= \frac{\partial \hat G}{\partial \mu}(x,0)=a+\psi (x),
\]
with $\min (a+\psi (x))>0$, then there exists $\mu_0>0$ and $k_5>0$ such that $|\rho(G)-\rho(G_0)|\le k_5\mu^2$ for all $\mu\in (0,\mu_0]$.\label{cor:rho2}\end{corollary}

This follows from the fact that the conditions of Proposition~\ref{prop:rhoclose} hold and from \eqref{eq:rhoclose} $T_0$ is the same for both $G$ and $\hat G$ so the difference in the rotation number is at order $\mu^2$.

\begin{corollary}If $G$ is in class $\mathcal{R}_{0,1}$ and $\min (a+\psi (x))>0$ then $\rho(G)$ is differentiable at $\mu =0$  and the derivative is $T_0^{-1}$. \label{cor:0diff}\end{corollary}

\emph{Proof:}

Clearly if $\mu>0$ \eqref{eq:rhoclose} implies that the rotation number is differentiable from above and the derivative is $T_0^{-1}$. The extension to $\mu <0$ follows the same argument as in \cite{GMM2025}. If $\mu<0$ then 
\[
G(x,\mu )=x-|\mu |(a+\psi(x))+\mu^2g(x,-\mu ).S
\]
But $G^{-1}(x, \mu)=x-\mu (a+\psi (x))+\mu^2\hat g(x,\mu )$ for some bounded periodic function $\hat g$, i.e.
\[
|G(x,-\mu )-G^{-1}(x,\mu )|=\mu^2(g(x,-\mu )-\hat g(x,\mu)).
\]
Now, $\rho (G^{-1})=-\rho (G)$, and Corollary~\ref{cor:rho2} imply that there exists $k_6>0$ such that $|\rho (G(x,-\mu )+\rho (G(x,\mu))|<k_6\mu^2$, so $|\,\rho (x,-\mu)+T_0^{-1}|\mu |\,|<k_6\mu^2$ or
$|\rho (x,-\mu)-T_0^{-1}\mu |<k_6\mu^2$.
\newline\rightline{$\square$}


\section{Proof of Theorem~\ref{thm:main}}\label{sect:proofmain}
By Proposition~\ref{prop:FS} the $q^{th}$ iterate of a lift of the map takes the form \eqref{eq:FqFS} which satisfies the hypotheses of Proposition~\ref{prop:rhoclose} once an integer is subtracted. Thus $G(x,\mu )=F^q(x,\mu )-p$ is a family of maps in class  $\mathcal{R}_{0,1}$ with
\[
\frac{\partial G}{\partial \mu}(x,0)=qa+q\Psi (x),
\]
with $\Psi$ defined by \eqref{eq:dFqmu}. The family $G$ are lifts of the $q^{th}$ iterate of the underlying family of circle diffeomorphisms.  

If $\min (a+\Psi)>0$ then by Corollary~\ref{cor:0diff}, $\rho (G (x,\mu))$ is differentiable as a function of $\mu$ at $\mu =0$ and  
\[
\frac{\partial \rho (G)}{\partial \mu}\big|_{\mu=0}=T_G^{-1}, \quad T_G=\frac{1}{q}\int_0^1\frac{1}{a+\Psi (x)}dx=\frac{1}{q}T_0.
\]
From elementary manipulation, $\rho (F^q)=p+\rho (G)$ and $\rho (F)=\frac{1}{q}\rho (F^q)$. Hence $\rho (F)=\rho (\mu)$ is differentiable and the derivative is given by \eqref{eq:main} as required for the first part of the theorem.

Now suppose that $\min (a+\Psi (x))<0<\max (a+\Psi (x))$. Let $y_1$ and $y_2$ satisfy $|y_1-y_2|<1$ and $a+\Psi(y_1)=\min (a+\Psi (x))$ and $a+\Psi(y_2)=\min (a+\Psi (x))$. Then $a+\Psi(y_1)<0$ and $a+\Psi(y_2)>0$ and hence there exists $x_1$ with $a+\Psi (x_1)=0$. By continuity there is a (possibly smaller) open neighbourhood of $\mu =0$ such that in the notation of \eqref{eq:Fq} $a+\Psi(y_1)+\mu g_q(y_1,\mu )<0$ and $a+\Psi(y_2)+\mu g_q(y_2,\mu )>0$ and hence there is a fixed point of $G$, which is periodic point of period $q$ for $F$, between $y_1$ and $y_2$ for all $\mu$ in that open neighbourhood. So the rotation number is constant and equal to $\frac{p}{q}$ on this neighbourhood, and is differentiable with derivative zero at $\mu=0$.
\newline\rightline{$\square$}

\medskip
Note that in the case with $c_{nq}=0$ for all $n\in \mathbb{Z}\backslash\{0\}$, $\Psi_q\equiv 0$, the constant $T_0=\frac{1}{qa}$ for $F^q$ and so we recover the same value for the derivative for the rotation number of $F$ (i.e. $\frac{1}{q}T_0^{-1}=a$) as would be obtained using the formula of Theorem~\ref{thm:H1977}(b) due to Brunovsky \cite{Brunovsky1974} for the case of irrational rotation number.

\section{The Arnold circle map revisited}\label{sect:morearnold}

In Sec.~\ref{sect:Arnold} only the simplest cases of scaling were considered, those cases for which there are no resonant terms in the Fourier series of the map and hence $\Psi(x)\equiv 0$ and the calculation of the derivative is trivial, see Corollary~\ref{cor:Arnold}. In this section we will look at the more complicated case with $(p,q)=(0,1)$, for which $\Psi(x )=\psi (x)=\beta \sin (2\pi x)$ and a modified Arnold map which includes higher order harmonics. 

\begin{corollary}
Consider a family of Arnold maps \eqref{eq:Arn00} with coefficients $(\alpha (\mu ),\beta (\mu ))$ such that $(\alpha (0),\beta (0))=(0,0)$. If $|\alpha '(0)|>|\beta '(0)|$, where primes denote differentiation with respect to $\mu$,  then the rotation number $\rho (\mu )$ is differentiable at $\mu =0$ and
\begin{equation}
\frac{d\rho}{d\mu}(0)={\rm sign}(\alpha '(0))\sqrt{|\alpha '(0)|^2-|\beta '(0)|^2}.
\label{eq:Arndiff0}\end{equation} 
\label{cor:0diffA} \end{corollary}

\emph{Proof:}

It is an exercise to check that the conditions of Proposition~\ref{prop:rhoclose} hold if $|\alpha '(0)|>|\beta '(0)|$ and then \eqref{eq:Arndiff0} follows from direct evaluation of $T_0$ using the standard definite integral
\begin{equation}
\int_0^1\frac{1}{a+b\sin (2\pi x)}dx=\frac{1}{\sqrt{a^2-b^2}}, \quad {\rm if}~a>b>0.
\label{eq:stanint}\end{equation}

\rightline{$\square$}

A numerical study of \eqref{eq:Arn00} on two paths through $(0,0)$, $(\alpha (\mu ),\beta (\mu ))=(5\mu ,\mu )$ and $(\alpha ,\beta )=(5\mu ,4\mu)$, give approximate straight lines for $\rho (\mu )$ with slopes $4.9$ (compared to the predicted $\sqrt{24}\approx 4.898$) and $3$ (the predicted value) respectively. (Rotation numbers were calculated using $10^5$ iterations of the map at 100 points between  $\mu=0$ and $\mu=0.1$. The slopes were computed using the least squares subroutine {\tt reglin} in Scilab \cite{Scilab}.) 

To investigate examples with higher order harmonics we will consider the modified Arnold map
\begin{equation}
G(x,\alpha ,\beta ,\gamma )=x+\alpha+\beta\sin (2\pi x) +\gamma\cos (4\pi x).
\label{eq:modArnold}\end{equation}

Theorem~\ref{thm:main} describes the scaling of the rotation number near $(\alpha,\eta,\gamma )=(\frac{p}{q},0,0)$ so we will make the parametrization explicit by writing
\begin{equation}
\alpha (\mu )=\frac{p}{q}+a\mu, \quad \beta (\mu )=b\mu, \quad \gamma (\mu) =c\mu,
\label{eq:newpar}\end{equation}
with $\mu$ passing through zero. In more general terms $\alpha'(0)=a$, $\beta'(0)=b$ and $\gamma'(0)=c$. 

We will consider two further cases, showing the different ways in which the terms interact.

If $(p,q)=(1,2)$ then $\Psi (x)=c\cos (4\pi x)$ and provided $a>|c|$ Theorem~\ref{thm:main} holds and using \eqref{eq:stanint}
\begin{equation}
\rho'(0)=\left(\int_0^1\frac{1}{a+c\cos (4\pi x)}dx\right)^{-1}=\sqrt{a^2-c^2}.
\label{eq:1exact}\end{equation}
Note that the derivative is independent of $b$. We will return to numerical confirmation at the end of this section. 

If  $(p,q)=(0,1)$ then $\Psi (x)=b\sin (2\pi x) +c\cos (4\pi x)$ and the transversality condition \eqref{eq:genericity} is clearly 
satisfied if $a>|b|+|c|$. To make the derivative of $\rho$ easier to evaluate we will consider a special choice of $a$, $b$ and $c$ which allows the integrals to be solved exactly. For $u\in (-\frac{1}{2},\frac{1}{2})$ define
\begin{equation}
a=1+\frac{1}{2}u^2, \quad b=2u, \quad c=-\frac{1}{2}u^2.
\label{eq:udef}\end{equation}

\begin{lemma}Consider \eqref{eq:modArnold} with parameters defined by \eqref{eq:newpar} and \eqref{eq:udef}. If $(p,q)=(0,1)$ then for each fixed $u\in (-\frac{1}{2},\frac{1}{2})$ the rotation number $\rho (\mu )$ is differentiable at $\mu=0$ and
\begin{equation}
\rho'(0)=(1-u^2)^{\frac{3}{2}}.
\label{eq:2exact}\end{equation}
\label{lem:01mod}\end{lemma}

\emph{Proof:}
If \eqref{eq:udef} holds and $u\in (-\frac{1}{2},\frac{1}{2})$ then $a>|b|+|c|$, so \eqref{eq:genericity} holds and Theorem~\ref{thm:main} is applicable with
\[\begin{array}{rl}
a+\Psi (x) &= 1+\frac{1}{2}u^2+2u\sin (2\pi x)-\frac{1}{2}u^2\cos (4\pi x)\\ &=(1+u\sin (2\pi x))^2.
\end{array}
\]
Thus for each $u\in (-\frac{1}{2},\frac{1}{2})$, $\rho (\mu )$ is differentiable at $\mu =0$ and
\[
\rho'(0)=\left( \int_0^1\frac{1}{(1+u\sin (2\pi x))^2}dx\right)^{-1},
\]
which equals \eqref{eq:2exact} as this is a standard integral in standard form.
\newline
\rightline{$\square$}

A small number of numerical experiments have been conducted to bring out these features, The results are summarized in 
Table~\ref{table}. Each row of the Table describes the parameters chosen for the study of \eqref{eq:modArnold} as a function of $\mu$ with parameters either given in terms of a triple $(a,b,c)$ representing the derivatives of $\alpha$, $\beta$ and $\gamma$ at $\mu =0$ via \eqref{eq:newpar}, or the value of $u$ from \eqref{eq:udef}.  If $c=0$ then the results are relevant to the original Arnold map \eqref{eq:Arn00}. The only anomaly in Table~\ref{table} is the value of the numerically determined derivative on the third line, with $(a,b,c)=(5,4,0.5)$. Running the same experiment using 20 points evenly spaced in $(0,0.001]$ instead of $(0,0.01]$ gives the estimate $\rho'(0)\approx 4.977$, so the interval on which the derivative is calculated needs to be smaller to obtain an accurate value. We assume that this is because of the large relative size of $b$: the derivative is independent of $b$ but higher order terms, and hence the neighbourhood on which the linear approximation is a good approximation, will be sensitive to the size of $b$. This effect is particularly marked at larger values of $q$.   

\begin{table}
\centering
\begin{tabular}{ |c|c|c|c|c| } 
 \hline
 $(p,q)$ & $(a,b,c)$ or $u$ & $\rho'(0)$ numerical & $\rho'(0)$ theoretical & equation\\ 
 \hline\hline
$(0,1)$ & $(5,4,0)$ & 3.003 & 3 & \eqref{eq:Arndiff0}\\ 
$(0,1)$ & $(5,1,0)$ & 3.003 & 3 & \eqref{eq:Arndiff0}\\ 
 & & & & \\
$(1,2)$ & $(5,4,0.5)$ & 5.007 & 4.975  & \eqref{eq:1exact}\\
$(1,2)$ & $(5,-1,0.5)$ & 4.978 & 4.975  & \eqref{eq:1exact}\\
$(1,2)$ & $(5,1,3)$ & 4.026 & 4  & \eqref{eq:1exact}\\
 & & & & \\
$(0,1)$ & 0.2 & 0.940 & 0.941 & \eqref{eq:2exact}\\ 
$(0,1)$ & $-0.2$ & 0.940 & 0.941 &  \eqref{eq:2exact}\\ 
$(0,1)$ & 0.4 & 0.770 & 0.770 &  \eqref{eq:2exact}\\ 
$(0,1)$ & $-0.4$ & 0.770 & 0.770  &  \eqref{eq:2exact}\\

 \hline
\end{tabular}
\caption{Results of numerical computation of $\rho'(0)$ for \eqref{eq:modArnold} compared with the theoretical prediction. In all cases the rotation number was calculated at 20 values of $\mu$ evenly spaced in $(0,0.01]$ using $10^5$ iterations of the map at each parameter value. The slope was calculated using the {\tt reglin} linear regression subroutine in Scilab \cite{Scilab}. Numbers are entered to three decimal places and if a scalar ($u$) is entered in the second column the parameters are defined via \eqref{eq:udef}. The first two lines are included in the text below the proof of Corollary~\ref{cor:0diffA} and are included here for completeness.}
\label{table}
\end{table}


\section{Piecewise linear maps}\label{sect:PWL}
Piecewise linear maps have many interesting features\cite{Avrutin2019}. Piecewise linear circle homeomorphisms can provide examples that are more amenable to analysis than smooth models, and they were studied by Herman \cite{H1979} amongst others. More recently the general piecewise smooth cases have been considered in some detail for irrational rotation numbers, e.g. \cite{Dzh2018, Khanin2003, Liousse2005}. The structure of mode locked regions and applications have also attracted interest \cite{Campbell1996, Carretero1997,Gaivao2025, JKS2018, Simpson2018}. The problem of differentiability of the rotation number is addressed in \cite{GMM2025} along with other properties, and the approach here can also be adapted to this setting.
    
Families of lifts of piecewise linear circle maps can be defined in a neighbourhood ${\mathcal N}$ of $\mu =0$ via $N$ pairs of $C^2$ functions $(a_k(\mu ), b_k(\mu ))$. The functions $b_k$ determine the break points of the map. These are the places at which the slope of the map changes and we chose the labelling so that
$b_1 (\mu )<b_2(\mu )<\dots <b_N(\mu )<b_{N+1}(\mu )=1+b_1(\mu)$ (extending the labelling in the obvious way). The functions $a_k(\mu )$ are the values of the map at the corresponding break point with $a_{N+1}=a_1+1$. The corresponding slopes $s_k(\mu )$ are defined by
\begin{equation}
    s_k(\mu )=\frac{a_{k+1}(\mu)-a_k(\mu )}{b_{k+1}(\mu) -b_k(\mu )},
\label{eq:slopes}
\end{equation}
$k=1, \dots ,N$. The lifts are therefore defined as
\begin{equation}
F(x,\mu )=a_i(\mu )+s_i(\mu )(x-b_i(\mu)), \quad b_i\le x< b_{i+1},   
\label{eq:pwleq}\end{equation}
$i=1,\dots , N$, extended to the real line via $F(x+1,\mu )=F(x,\mu )+1$. This defines a homeomorphism with break points at $b_k$ provided
\begin{equation}
 s_i(\mu )>0 \quad {\rm and}\quad s_i(\mu )\ne s_{i+1}(\mu ),
\label{eq:pwlhom}
\end{equation}
for all $\mu\in{\mathcal N}$ and $i=1,\dots , N$. Such maps are described in \cite{GMM2025} in which it is shown that there are conditions on the orbits of the break points which ensure that $F^q$ is a rigid rotation through $p/q$ if $\mu =0$. Under these circumstances there is a neighbourhood of $\mu =0$ such that $F^q$ has $N_q$ break points $\gamma_i$ and there are constants $A_i$ and $B_i$, and a bounded periodic function $g$ as in Sec.~\ref{sect:FS} such that
\begin{equation}
F^q(x,\mu )=x+p+\mu (A_i+B_i(x-\gamma_i(\mu)) +\mu^2g(x,\mu ), 
\label{eq:pwleqq}\end{equation} 
$\gamma_i\le x< \gamma_{i+1}$, $i=1, \dots, N_q$. 
To avoid mode locking we impose the monotonicity condition in parameter space: 
\begin{equation}
\mu_1<\mu_2\ \ {\rm implies \ that}\ \ F(x,\mu_2)-F(x,\mu_1) >0,\ \ {\rm for\ all}\ x\in\mathbb{R}. 
\label{eq:nomodelock}\end{equation}
Condition \eqref{eq:nomodelock} can be written explicitly in terms of the movement of the breakpoints of $F$\cite{GMM2025}. The approach to this situation in \cite{GMM2025} is from first principles. They prove the following result.

\begin{theorem}\cite{GMM2025} Suppose that at $\mu=0$ there exists $\frac{p}{q}$ with $p$, $q$ coprime such that a family of piecewise linear circle homeomorphisms \eqref{eq:pwleq} has $F^q(x,0)=x +p$. If $F^q$ is written as \eqref{eq:pwleqq} and \eqref{eq:nomodelock} holds then $\rho (\mu)$ is differentiable at $\mu =0$ and
\begin{equation}
\rho'(0)=\frac{1}{q\Sigma_{i=1}^{N_q} T_i},
\label{eq:pwldiff}\end{equation}
where
\begin{equation}
T_i=\begin{cases}\frac{(\gamma_{i+1}-\gamma_i)}{A_i}& {\rm if}~~ B_i=0, \\
\frac{1}{B_i}\log \left(1+\frac{B_i(\gamma_{i+1}-\gamma_i)}{A_i}\right)& {\rm if}~~ B_i\ne 0.\end{cases}
\label{eq:pwlTi}
\end{equation}
\label{thm:GMM}
\end{theorem}

The methods described in earlier sections of this paper provide an alternative proof of Theorem~\ref{thm:GMM} which we sketch here.

We consider the map $F^q(x)-p$ from \eqref{eq:pwleqq} in a neighbourhood of $\mu =0$. This is approximated using Euler's method by the linear ordinary differentials
\begin{equation}
\frac{dX}{dt}=A_i+B_i(X-\gamma_i), \quad \gamma_i\le X\le \gamma_{i+1}.
\label{eq:pwlode}
\end{equation}
By simple integration the time taken from $\gamma_i$ to $\gamma_{i+1}$ is therefore $T_i$ given by \eqref{eq:pwlTi}. The total time to move from $0$ to $1$ for $X$ is therefore $\sum T_i$ which plays the role of $T_0$ in Sec.~\ref{sect:int}. The methods of Sec.~\ref{sect:int} and ec.~\ref{sect:proofmain} can be used with the inequality
\[
e^{nh}-(1+h)^n\le nh^2e^{nh} 
\]
to control the errors and show that $\rho (F^q-p)=\frac{\mu}{\sum T_i}+O(\mu^2)$. Hence 
\[
\rho (\mu)=\frac{p}{q}+\frac{\mu}{q\sum T_i}+O(\mu^2),
\]  
completing the proof of Theorem~\ref{thm:GMM}.

A corollary of the differentiability of the rotation number in the cases we have considered is that there is no mode locking at that rational rotation number (cf. \cite{GMM2025}). We conjecture that this is linked to the existence of `pinching' in the locus of 
mode locked regions in two-parameter families of maps. This refers to situations in which the boundaries of mode locked regions intersect transversely, so that the interior of a mode locked region with a given rational rotation number is a series of non-intersecting regions whose closure is connected. This is typical in two parameter families of piecewise smooth maps with two break points \cite{Campbell1996}, but not three \cite{Boyland2025}.

\section*{Acknowledgements} I am grateful to Phil Boyland (Florida), Siyuan Ma and James Montaldi (Manchester), and David Simpson (Massey, Palmerston North) for helpful conversations.

\nocite{*}

\end{document}